\documentclass[11pt]{article}%
\usepackage{amsmath}
\usepackage{amsfonts}
\usepackage{amssymb}
\usepackage{color}
\usepackage{graphicx}
\providecommand{\U}[1]{\protect\rule{.1in}{.1in}}
\providecommand{\U}[1]{\protect\rule{.1in}{.1in}}
\newtheorem{theorem}{Theorem}[section]

\newtheorem{definition}{Definition}[section]

\newtheorem{remark}{Remark}[section]

\renewcommand{\theequation}{\arabic{section}.\arabic{equation}}
\newenvironment{m_proof}[1][Proof]{\noindent\textit{#1.} }{\ \hfill\rule{0.3em}{0.5em}}
\def\Pr{{\rm Pr}}
\def\E{{\rm E}}

\begin{document}

\voffset-1cm

\title{\bf On Hodges' Superefficiency and Merits of Oracle Property in Model Selection\thanks{This work was partially supported by NSFC (grant No. 71771089) and the 111 Project (grant No. B14019).}}
\author{ Xianyi Wu and Xian Zhou\\  East China Normal University, Shanghai,  China\\
Macquarie University, Sydney, Australia.
}
\date{}

\maketitle
\begin{abstract}

The oracle property of model selection procedures has attracted a large volume of favorable publications in the literature, but also faced criticisms of being ineffective and misleading in applications. Such criticisms, however, have appeared to be largely ignored by the majority of the popular statistical literature, despite their serious impact.  In this paper, we present a new type of Hodges' estimators that can easily produce model selection procedures with the oracle and some other desired properties, but can be readily seen to perform poorly in parts of the  parameter spaces that are fixed and independent of sample sizes. Consequently, the merits of the oracle property for model selection as extensively advocated in the literature are questionable and possibly overstated. In particular, because the mathematics employed in this paper are at an elementary level, this finding leads to new discoveries on the merits of the oracle property and exposes some overlooked crucial facts on model selection procedures.  

\bigskip
\noindent\textbf{Keywords}: Hodges' estimator; model selection; oracle property; penalized maximum likelihood/least squares;  superefficiency.
\end{abstract}

%\newpage

\section{Introduction}

Model selection is undoubtedly an extensively employed technique in data analysis and has attracted a great deal of research interests in the literature. It has become increasingly popular and ubiquitous partly due to rapid advance in computational power. A large volume of literature on model selection has been published, and widely used statistical software packages contain more or less routines for model selection. In particular, almost all textbooks on linear regressions (including those at undergraduate level) dedicate a whole chapter or a few sections to variable (model) selection.

Theoretically, there are two purposes for model selection in general: model identification and inference. The latter is often referred to as post-selection inference in practical data analyses (Berk et al, 2010, 2013).  Typically in post-selection inference, model identification is an intermediate step in data analysis, after that the analysts then perform further statistical analyses based on the selected model, pretending it is the ``true'' model.

This amounts to a two-step procedure in data analysis (Berk et. al. 2010, 2013). 
Many textbooks on linear regression oriented in practical data analysis taught students how to use well-known techniques such as AIC (Akaike's information criterion), BIC (Bayesian information criterion), adjusted R-squares, Mallows' C$_p$ and so on to select ``good'' or ``best'' models, then draw statistical inference based on the selected models. Popular softwares, such as SAS and R, generally select a model by some criterion (such as AIC/BIC) and output the estimated coefficients of the selected explanatory variables as well as their corresponding $p$-values computed as if the selected models truly represent the nature. Berk et al. (2010) provides an exemplified but limited list of remarkable research works of this type. More and interesting discussions on how model selection methods perform for model identification and inference can be found in Yang (2005, 2007).

In this paper, we focus on the parameter estimation for model selection based on an overall model with a fixed parameter space of dimension $p$ (less than the sample size $n$). The models are selected by setting some of its parameters to zero. Thus selecting a model corresponds to estimating some parameters by zero. Statistical inference is drawn from the nonzero estimates of the other parameters on the selected model. In this sense, the parameter estimation we consider covers both model identification and inference.

More recently, based on the idea of penalized maximum likelihood/least squares inherited from AIC and BIC, some researchers found a smart way to integrate this two-step analysis into a single-step procedure. It associates some cleverly designed penalties to the likelihood function (or squared errors in regression models) so that, by maximizing the penalized likelihood function (or minimizing the penalized squared errors in regression models), a part of parameters are estimated by zero and the others by nonzero quantities. This enables model identification and inference to be carried out together -- referred to as a one-step procedure for model selection.  

This new approach began with the famous LASSO algorithm (Least Absolute Selection and Shrinkage Operator) proposed by Tibshirani (1996). Generally, let $l(Y, \theta)$ denote the log-likelihood of the parameters $\theta=(\theta_1,\theta_2,\dots,\theta_p)'$. Then the estimators of $\theta_i$s are obtained by maximizing the penalized log-likelihood 
$
PL(\theta)=l(Y, \theta)+\sum_{i=1}^pf_i(\theta_i, \lambda_n),
$
where $Y$ is the sample of size $n$, $f_i(\theta_i, \lambda_n)$ the penalty associated with parameter $\theta_i$, and $\lambda_n$ a known tuning parameter. In parallel, the estimators for a linear regression model 
$Y=X_{n\times p}\theta+\varepsilon$  are obtained by minimizing penalized squares 
$
PS(\theta)=(Y-X\theta)'(Y-X\theta)+2\sum_{i=1}^pf_i(\theta_i, \lambda_n).
$
For particularly designed penalty $f_i$, maximizing $PL(\theta)$ produces such estimators that automatically estimate some $\theta_i$ by zero, so as to simultaneously select model and estimate the  parameters of the selected model. 
Different designs of the penalties generate different selection-estimation algorithms, such as
LASSO (Tibshirani, 1996), adaptive LASSO (Zou, 2006), hard thresholding estimators, soft thresholding estimators, bridge regression estimators (Frank and Friedman, 1993), SCAD (smoothly clipped absolute deviation penalty, Fan and Li, 2001), Elastic net method (Zou and Hastie, 2005) and MCP (minimax concave penalty, Zhang, 2010). 

The first theoretical justification of the penalized maximum likelihood method
appeared in Fan and Li (2001), who proved an appealing property of their SCAD estimators, referred to as the {\it oracle property}, as defined below. 

\begin{definition}\label{oracle property} Let
$b(\theta)=\{i: \theta_i\neq 0,i=1,2,\dots,p\}\hbox{ and }\bar b(\theta)=\{1,2,\dots,p\}-b(\theta)$
and rearrange the parameter vector as 
$\theta=(\theta_{b(\theta)}',\theta_{\bar b(\theta)}')'$ without loss of generality.
An estimator $\hat\theta_n$ is said to possess the oracle property or an oracle estimator if
\begin{enumerate}
\item[(1)] $\lim_{n\rightarrow\infty}\limits\Pr_\theta(\hat\theta_{n,\bar b(\theta)}=0)=1$; and
\item[(2)] $\sqrt n (\hat\theta_{n,b(\theta)}-\theta_{b(\theta)})\overset{d}{\rightarrow} N(0, {\cal F}_{b(\theta)}^{-1})$, where $\overset{d}{\rightarrow}$ indicates convergence in distribution  and ${\cal F}_{b(\theta)}$ is the Fisher's information matrix of $\theta_{b(\theta)}$ knowing $\theta_{\bar b(\theta)}=0$.
\end{enumerate}
 \end{definition}
The oracle property defined above is described in limit sense as the sample size tends to infinity. It states that an oracle estimator is asymptotically equivalent to the ideal estimator with the active (nonzero) parameters only. This has extensively been promoted as a justification to using the single-step procedure  SCAD estimation. The follow-ups of Fan and Li (2001) can be summarized in two aspects:
\begin{enumerate}
\item[(1)] Many researchers  have focused on finding model selection procedures for linear regression so as to produce oracle estimators that perform better in their finite sample simulations. Significant works in this line include adaptive LASSO by Zou (2006), adaptive group LASSO by Wang and Leng (2008), elastic net method by Zou and Hastie (2005), adaptive elastic net method by Zou and Zhang (2009), MCP by Zhang (2010), Orthogonalising EM algorithm with nonconvex penalties by Xiong et al. (2017) and other variants by, e.g., Leng et al. (2014), Gefang (2014) and Wang et al. (2011), among a vast number of others.  
\item[(2)]  More papers have aimed at extending the  ideas growing up in linear regression to other models so as to derive estimators with the oracle property.  Important examples include generalized linear models (Friedman, et al., 2010; van der Geer, 2008), Cox regression models (Fan and Li, 2002; Zhang and Lu, 2007), network exploration (Fan et al, 2009), additive models (Horowitz et al., 2006; Huang et al., 2010),  panel data models (Fan and Li, 2004 and Kock, 2013) and propensity score models (Brookhart et al., 2006), among others.
\end{enumerate}
What is worth special mention is the explosion of the literature in model selection with the oracle property under a huge number of statistical models, which have led to thousands of citations to many of the papers mentioned above by following papers and books. 
 
Such a popular property, however, is not universally accepted without criticisms. The Vienna school, led by H. Leeb and B. M. P\"{o}tscher, argued that the oracle property is ``too good to be true'' and seriously questioned the advertised merits of oracle model selection procedures in a series of papers, including  Leeb and P\"{o}tscher (2008a, 2008b),  P\"otscher,  (2009), P\"{o}tscher and Leeb (2009), P\"{o}tscher and Schneider (2009) etc. and the references therein. They considered oracle estimators as a return  of Hodges' estimators (Le Cam, 1953), which provided a typical counterexample to a conjecture by R. A. Fisher that the MLEs have minimum asymptotic variance (referred to as asymptotic efficient), and are known to perform poorly at some parameter values in finite sample size. Their arguments are summarized below.
 
\begin{enumerate}
\item[(1)] Theoretically,  the overall risk  of a sparse estimator can be unbounded as $n\rightarrow \infty$.
\item[(2)] Numerous Monte Carlo simulations were conducted to demonstrate that,  in finite sample size, a number of oracle estimators perform quite poorly when the parameters take values in a vicinity of zero.
\end{enumerate} 
In particular, by replicating and extending Monte Carlo simulations of the performance of the SCAD estimator in Example 4.1 of Fan and Li (2001), Leeb and P\"otscher (2008b) demonstrated that this estimator, when tuned to enjoy sparsity and oracle properties, can perform poorly in finite samples. Even if not tuned to sparsity, the SCAD estimator performs worse than the least squares estimator in parts of the parameter space.
It is interesting to note that the simulation study in Fan and Li (2001) was conducted only at some points that happened to avoid the parts of the parameter space examined by Leeb and P\"otscher (2008b). This phenomenon highlights the fact that simulations may produce results to support either side of a controversy, thus are unreliable to judge the goodness of an estimator. This is generally recognized in the statistical community but often overlooked (more details will be discussed in ``The power of simulations'' in Section \ref{concluding}). Consequently, Leeb and P\"otscher (2008b) argued that the oracle property is highly misleading and cannot be relied upon to justify an estimator. 

It is also worth to note that a procedure with oracle property is consistent for model selection. While consistency is an important and desirable property in limiting sense for large-samples, it is not sufficient to justify the superiority of a procedure for finite (fixed) sample sizes. As indicated in Yang (2005, 2007), BIC type  procedures are consistent and AIC type procedures are inconsistent but optimal in minimax rate of convergence. Yang (2005, 2007) also argued that the strength of AIC and BIC could not be shared. More introduction on optimal model selection in minimax rate of convergence can be found in the recent work of Wang et al. (2014).

In this paper, we attempt to address this controversy by revisiting the issue of Hodges' superefficiency and theoretically analyzing the performance of oracle estimators in a class of generalized Hodges' estimators without relying on numerical simulations. As simulations are subject to large variations and may produce different results in either side of a controversy, they are not capable of resolving the controversy convincingly. The theoretical analysis in this paper can avoid such drawbacks of the simulation approach and thus provide more convincing conclusions.
In particular, because the mathematics employed in this paper are at a quite elementary level, this finding leads to some new discoveries on the merits of the oracle property and exposes some overlooked crucial facts on model selection procedures. In addition, some significant but often ignored facts on asymptotic properties are also emphasized to warn the pitfall of justifying the merits of statistical procedures based on asymptotic measures that the oracle property relies on.

Our arguments proceed as follows.
\begin{enumerate}
\item[(1)] Generally, asymptotic bias and variance do not reflect their versions in finite sample size. Thus the asymptotic superefficiency and the oracle property do not necessarily lead to good performance of the estimators in any fixed sample size.
\item[(2)] 
We introduce a new type of Hodges' estimators, referred to as {\it oracle Hodges' estimators}, which can reduce the asymptotic variance of any given estimator over certain multi-dimensional subsets of the parameter space and provide an easy way to obtain oracle model selection procedures.
\item[(3)] 
By analyzing the performance of oracle Hodges' estimators in finite sample size, we are able to theoretically demonstrate that such estimators perform poorly at some (true) values of the parameters despite possessing the oracle property. 
\item[(4)]
The performance of oracle Hodges' estimators indicates that the oracle property alone does not justify the use of oracle model selection procedures, and hence the merits of the oracle property as advocated in the literature are questionable. In particular there has been no rigorous theory, in any case where the MLE/LSE are available, to prove that oracle estimators are better than MLE/LSE under certain commonly accepted criteria (such as smaller mean squared error) with fixed sample size. 
\item[(5)] 
Our results provide theoretical justification to support the view of Leeb and P\"otscher (2008b) regarding the oracle model selection approach, and clarified that the oracle property is not a simple return of classical Hodges' superefficiency as suggested by Leeb and P\"otscher. Instead, the newly defined oracle Hodges' estimators, rather than their classical versions, have the same asymptotic behavior as the estimators from oracle model selection procedures.
\end{enumerate}

In Section \ref{section hodges} next, the new oracle Hodges' estimators are defined and their asymptotic behaviors (superefficiency) are discussed. Section \ref{model selection} introduces model selection methods derived from Hodges' estimators and discusses their connections to and differences from penalized maximum likelihood or least squares estimators. The performance of oracle Hodges' estimators is theoretically analyzed in Section \ref{performance} after a brief discussion on the connection of asymptotic biases and variances to their finite sample size versions. Some concluding remarks are discussed in Section \ref{concluding}.

%%%%%%%%%%%%%%%%%%%%%%%%%%%%%%%%%%%%%%%%%%%%%%%%%%%%%%%%%%%%%%%

\section{Asymptotic efficiency and Hodges' estimators}\label{section hodges}
\setcounter{equation}{0}
The concept of asymptotic efficiency was introduced by Sir R. A. Fisher with the attempt to justify the goodness of  MLEs and has great impacts on statistical inference in large samples. 
For a $d$-dimensional parameter $\theta=(\theta_1,\theta_2,\dots,\theta_d)'$, let $\hat\theta_n=(\hat\theta_{n1},\hat\theta_{n2},\dots,\hat\theta_{nd})'$ be any sequence of its estimators  such that $r_n(\hat\theta_n-\theta)\overset{d}{\rightarrow} Z$ for some sequence of scalars $r_n\rightarrow\infty$ and a  $d$-dimensional random vector $Z$ with mean zero and variance-covariance matrix $V^{-1}$ for some matrix $V=(v_{ij})_{d\times d}$. 
A sequence of estimators $\{\mathbf{\hat{\theta}}_{n}\}$ of $\mathbf{\theta}$ such that 
$\sqrt{n}(\mathbf{\hat{\theta}}_{n}-\mathbf{\theta)}\overset{d}\rightarrow N(0,\Sigma(\mathbf{\theta}))$, whatever is the true value of $\theta$, is said to be {\it asymptotically efficient} (or {\it asymptotically optimal} in some literature) if
$\Sigma(\mathbf{\theta})=\mathfrak{F}^{-1}(\theta)$, where $\mathfrak{F}(\theta)$ is the Fisher's information matrix of the distribution. Namely, $\hat\theta_n$ is efficient if it is asymptotically unbiased (in the sense that the limiting random variable $Z$ has zero mean)  and ``optimal'' (in a certain sense based on the asymptotic variance) in the class of asymptotically unbiased estimators with order $n^{-1/2}$.

The general results state that, under certain regularity conditions, a sequence of roots of likelihood equations
 is asymptotically efficient. There have been a huge number of research efforts aimed at seeking asymptotically efficient estimators. 
The idea of measuring efficiency by asymptotic variance, however, appears not as successful as its counterpart in finite sample size, and the lower bounds defined by $\mathfrak{F}^{-1}(\theta)$ are not sufficiently low in the class of asymptotically unbiased estimators. The example in Subsection \ref{Hodges' example} below, whose prototype was made widely known by Le Cam (1953) under the name Hodges' estimator, shows that, given any sequence of estimators $\hat\theta_n$ with $r_n(\hat\theta_n-\theta)\overset{d}{\rightarrow} Z$, one can construct another estimator sequence with asymptotic variance no more than that of $\hat\theta_n$ at any value of $\theta$, and strictly less at certain values of $\theta$. This phenomenon is referred to as (asymptotic) {\it superefficiency}. In spirit of this idea, we introduce variants of Hodges' estimators in Subsection \ref{Hodges' estimators in multiple dimensions}, which possess the oracle property (Definition \ref{oracle property}) and can easily serve the purpose of finding oracle model selection procedures.

Let $\|x\|=\sqrt{\sum_{i=1}^dx_i^2}$ denote the Euclidean norm of $x=(x_1,x_2,\dots,x_d)\in \mathbb{R}^d$, $d(x,A)=\inf_{y\in A}\|x-y\|$ the distance between a point $x$ and a subset $A$ and $d(A,B)=\inf_{x\in A,y\in B}\|x-y\|$  the distance between two subsets of $\mathbb{R}^d$ under the Euclidean norm.

%%%%%%%%%%%%%%%%%%%%%%%%%%%%%%%%%%%%%%%%%%%%%%%%%%%%%%%%%%%%%%%%%%%%%%%%

\subsection{Classical Hodges' example of superefficiency}\label{Hodges' example}

Let $\hat{\theta}_{n}$ be any sequence of estimators such that
$r_{n}(\hat{\theta}_{n}-\theta)\overset{d}{\rightarrow}Z$ for a sequence of deterministic scalars
$r_{n}\rightarrow\infty$, where the distribution of $Z$ may depend on $\theta$, 
whatever is the true value of $\theta$. 

\begin{definition}[Hodges' estimator]\label{Hodges' estimators_naive}
Let $\{a_n\}$ be a sequence of scalars and $c$ any fixed point in the parameter space $\Theta$. The Hodges' estimator of $\theta$ is defined by
\begin{equation}\label{Hodges_one_D}
\breve{\theta}_{n}(c)=cI(\|\hat{\theta}_{n}-c\|\leq a_{n})+\hat{\theta}_{n} I(\|\hat{\theta}_{n}-c\|>a_{n}).
\end{equation}
\end{definition}

\noindent This estimator was initiated by Jr. Hodges with $r_n=\sqrt{n}$, $c=0$ and $a_n=n^{-1/4}$ for one-dimensional $\theta$ (reported by Le Cam, 1953; see also Lehmann and Casella, 1998, p. 420, Example 2.5) and has been revisited many times by, e.g., Leeb and P\"otscher (2005, 2008a, b) when they discussed the implications of consistent model selections. 
While (\ref{Hodges_one_D}) may be slightly generalized to
\[
\breve{\theta}_{n}(c)=((1-\alpha)c+\alpha\hat\theta_n)I(|\hat{\theta}_{n}-c|\leq a_{n})
+\hat{\theta}_{n}I(|\hat{\theta}_{n}-c|>a_{n}),
\]
where $\alpha\in[0,1]$ (cf., e.g., van der Vaart, 1998 for the case with $r_n=\sqrt n$ and $c=0$), we here take $\alpha=0$ as this is sufficient for our purpose.
The asymptotic distribution of $\breve{\theta}_n(c)$ in Definition \ref{Hodges' estimators_naive} is given in the following theorem.

\begin{theorem}\label{asymp_1}
For any sequence $\{a_n\}$ such that 
$a_{n}=o(1)$ and $r_{n}a_n\rightarrow\infty$ as $n\rightarrow\infty$, we have
\[
r_{n}(\breve{\theta}_{n}(c)-\theta)\overset{d}{\rightarrow}ZI(\theta\not =c).
\]
\end{theorem}

It is clear that, in terms of asymptotic variances,  $\breve{\theta}_{n}(c)$ is no worse than $\hat{\theta
}_{n}$ at any $\theta$ and strictly better than $\hat\theta_n$ at $\theta=c$ because the asymptotic variance of $\breve\theta_n(c)$ is zero at $\theta=c$. This example revealed an interesting phenomenon that, in terms of the asymptotic variance, any estimate can be improved at an arbitrary but fixed point in the parameter space.

%%%%%%%%%%%%%%%%%%%%%%%%%%%%%%%%%%%%%%%%%%%%%%%%%%%%%%%%%%%%%%%%%%%%%%%%%%%%%%%%%%%%%%%%%%%%%%%

\subsection{Oracle Hodges' estimators}\label{Hodges' estimators in multiple dimensions}
We next introduce a few closely linked variants of Hodges' estimators, which differ from the classical version in Definition \ref{Hodges' estimators_naive}, but we have kept Hodges' name for the new estimators because they retain the feature of superefficiency. We will refer to the new type of Hodges' estimators as {\it oracle Hodges' estimators} because they possess the oracle property in Definition \ref{oracle property}, as will be shown via Theorems \ref{multi-dimension-super-efficiency} to \ref{LSE superefficiency} below. The definition of oracle Hodges' estimators is presented in four versions: general version, continuous/smoothing version, MLE version and LSE version, where the first one is fundamental, the second is a refinement, and the other two are special cases with origins $\hat\theta_n$ being MLE and LSE, respectively.
\vskip3mm

\noindent{\bf 1. General version.}

\noindent For any subset $b\subset \{1,2,\dots,d\}$ and its complement $\bar b=\{1,2,\dots,d\}-b$, rearrange the components of $\theta$ and $\hat\theta_n$ as $\theta'=(\theta_b',\theta_{\bar b}')$ and $\hat\theta_n'=(\hat\theta_{n,b}',\hat\theta_{n,{\bar b}}')$, respectively, such that $r_n(\hat\theta_n-\theta)\overset{d}{\rightarrow} Z$, where the current $Z=(Z_b',Z_{\bar b}')'$ is also a rearrangement of the original $Z$ in the same way as $\theta$. Clearly, the mean of $Z$ is still a zero vector but the covariance matrix changes to $V^{-1}$ with
\[
V=\left(\begin{array}{cc}V_{bb}&V_{b\bar{b}}\\V_{\bar{b}b}&V_{\bar{b}\bar{b}}\end{array}\right)
=\left(\begin{array}{cc}(v_{ij})_{i,j\in b}&(v_{ij})_{i\in b,j\in\bar b}\\
(v_{ij})_{i\in\bar b,j\in b} & (v_{ij})_{i,j\in\bar b}\end{array}\right).
\] 
It is also easy to see that, if $b\neq\emptyset$, then the marginal vector $\hat\theta_{n,b}$ has an asymptotic distribution given by $r_n(\hat\theta_{n,b}-\theta_{b})\overset{d}{\rightarrow} Z_b$ with mean zero and covariance matrix 
\begin{equation}\label{Delta_theta}
\Delta_b=(V_{bb}-V_{b\bar b}V_{\bar b\bar b}^{-1}V_{\bar b b})^{-1}.
\end{equation} 
Note that the distribution of $Z$, and hence the variance matrix $V^{-1}$, may be related to parameters $\theta$. The only requirement is that $V^{-1}$ is a continuous function of $\theta$. We use $\hat V^{-1}$ to denote any of consistent estimators of $V^{-1}$, e.g., obtained by substituting $\hat\theta_n$ for $\theta$, so that the symbols $\hat V_{bb}$, $\hat{V}_{b\bar b}$ and so on are self-explained.

Let  $c=(c_1,c_2,\dots,c_d)'$ be a known $d$-vector. For every nonempty and proper subset $b$ of $\{1,2,\dots,d\}$ (i.e., $\emptyset\ne b\ne \{1,2,\dots,d\}$), denote 
\begin{equation}\label{part_estimate}
\check\theta_{n, b}= \hat\theta_{n,b}+\hat{V}_{bb}^{-1}\hat{V}_{b\bar{b}}(\hat\theta_{n,\bar{b}}-c_{\bar b})
\quad\hbox{and}\quad\check\theta_n(b)=(\check\theta_{n,b}',c_{\bar b}')'
\end{equation}
with the convention $\check\theta_{n,\{1,2,\dots,d\}}=\check\theta_{n}(\{1,2,\dots,d\})=\hat\theta_{n}$. 
Moreover, we redefine $b(\theta)$ and $\bar{b}(\theta)$ by
\begin{equation}\label{b_theta_1}
b(\theta)=\{j:j\in\{1,2,\dots,d\},\theta_j\neq c_j\} \quad\hbox{and}\quad \bar{b}(\theta)=\{1,2,\dots,d\}-b(\theta).
 \end{equation} 
The following definition introduces a sequence $\tilde\theta_n (c)$ of oracle Hodges' estimators in multi-dimensional case derived from $\hat\theta_n$.

\begin{definition}[Oracle Hodges' estimators]\label{Hodges' estimators}
 {\rm Let $(a_{nj})=(a_{n1},\dots,a_{nd})$, $n=1,2,\dots$, denote a sequence of $d$-vectors with positive components.
 For every $n=1,2,\dots$, define two complementary random sets by
 \begin{equation}\label{b_n}
 {b}_n(c)=\{j: |\hat\theta_{nj}-c_{j}|>a_{nj}\}\quad\hbox{and}\quad
{\bar b}_n(c)=\{j: |\hat\theta_{nj}-c_{j}|\leq a_{nj}\}
 \end{equation}  
and the corresponding oracle Hodges' estimator by
\begin{equation}\label{Hodges' estimate-multi-dimension}
\tilde\theta_n(c)=\check\theta_n(b_n(c))=\left\{\begin{array}{ll}\hat\theta_n&\hbox{ if } b_n(c)=\{1,2,\dots,d\},\\(\check\theta_{n,b_n(c)}',c_{\bar b_n(c)}')'& \hbox{ if } b_n(c)\neq\emptyset,b_n(c)\neq\{1,2,\dots,d\},\\
c&\hbox{ if }b_n(c)=\emptyset,\end{array}\right.
\end{equation}
where $\check\theta_n(b_n(c))$ and $\check\theta_{n,b_n(c)}$ are obtained from the two equations in \eqref{part_estimate} by substituting $b_n(c)$ for $b$, and $b_n(c)$ is defined in \eqref{b_n}.
}
\end{definition}
For later reference, denote  
\begin{equation}\label{X_b}
\check Z_{b}=\left\{\begin{array}{ll}V_{{b}{b}}^{-1}(V_{{b}{b}}\quad V_{{b}\bar{b}})Z,&\hbox{ if } b\neq\emptyset,\\
0,&\hbox{ if } b=\emptyset.\end{array}\right.
\end{equation}
 The asymptotic properties of $\tilde\theta_n(c)$ defined in (\ref{Hodges' estimate-multi-dimension}) are provided in the next theorem.
 
\begin{theorem}\label{multi-dimension-super-efficiency}
{\rm If the sequence of $d$-vectors $\{(a_{nj})\}$ satisfies  
\begin{equation}\label{condition on a_n}
\max_{1\leq j\leq d} a_{nj}\rightarrow0\quad\hbox{ and }\quad
r_n\min_{1\leq j\leq d}a_{nj}\rightarrow\infty\hbox{ as }n\rightarrow\infty,
\end{equation}
then for any $b\subset\{1,2,\dots,d\}$,  $\theta_{\bar b}=c_{\bar b}$ implies 
$\lim_{n\rightarrow\infty}\limits\Pr(\tilde\theta_{n,\bar b}(c)=\theta_{\bar b})\rightarrow 1$ and
\begin{equation}\label{2.10}
r_n(\tilde\theta_n(c)-\theta)\overset{d}{\rightarrow} 
\left(\begin{array}{c}\check Z_{b(\theta)}\\ 0\end{array}\right) \hbox{ under }\Pr_\theta,
\end{equation} 
whatever is the true value of $\theta$, where $\check Z_{b(\theta)}$ is defined as in \eqref{X_b} with $b$ replaced by $b(\theta)$ in \eqref{b_theta_1}.}
\end{theorem}

\noindent{\bf 2. Continuous/smoothing version.}

\noindent
As a function of $\hat\theta_n$, the estimator $\tilde\theta_n(c)$ is not continuous  at any point  $\hat\theta_n=(\hat\theta_{n1},\hat\theta_{n2},\dots,\hat\theta_{nd})$ such that $|\hat\theta_{nj}-c_j|=a_{nj}$ for some $j\in\{1,2,\dots,d\}$. Some authors think of the continuity as an important property (see, e.g., Fan and Li, 2001) but others may disagree. If preferred, a continuous version of $\tilde\theta_n(c)$ can be easily achieved by the following procedure.

Let $(a_{nj}^{(1)}) $ and $(a_{nj}^{(2)})$ be two sequences of $d$-vectors both satisfying conditions \eqref{condition on a_n} on $(a_{nj})$ and $ a_{nj}^{(1)}< a_{nj}^{(2)}$, $j=1,2,\dots, d$. A possible choice is $a_{nj}^{(1)}=r_n^{-1/2}$ and $a_{nj}^{(2)}=2r_n^{-1/2}$, $j=1,2,\dots,d$. Define two sequences of oracle Hodges' estimators 
$\tilde{\theta}^{(1)}_n(c)$ and $\tilde{\theta}^{(2)}_n(c)$ by $(a_{nj}^{(1)}) $ and $(a_{nj}^{(2)})$ respectively as in Definition \ref{Hodges' estimators}. Let $f_1(x),\dots,f_d(x)$ be any $d$
continuous and increasing functions on $x\in\mathbb{R}^+$ such that $f_j(c_{nj}\pm a^{(1)}_{nj})=c_{nj}$ and 
$f_j(c_{nj}\pm a^{(2)}_{nj})=c_{nj}\pm a^{(2)}_{nj}$. Define
\begin{align}\label{continuous_Hodges}
\tilde\theta_{nj}(c;f)=\left\{\begin{array}{ll}c_{nj} & \hbox{ if }|\hat\theta_{nj}-c_{j}|\leq  a^{(1)}_{nj},\\
f(\hat\theta_{nj})& \hbox{ if }  a_{nj}^{(1)}\leq |\theta_{nj}-c_j |\leq  a_{nj}^{(2)},\\
\tilde{\theta}^{(2)}_{nj}(c)& \hbox{ otherwise.} \end{array}\right.
\end{align}
Then $\tilde\theta_{nj}(c;f)$ is a continuous version of $\tilde\theta_n(c)$ such that 
\[|\tilde\theta^{(2)}_{nj}(c)-c_{nj}|\leq |\tilde\theta_{nj}(c;f)-c_{nj}|<|\tilde\theta^{(1)}_{nj}(c)-c_{nj}|.\] 
These inequalities ensure the following result.
\begin{theorem}
The estimators $\tilde\theta_{n}(c;f)$ defined in \eqref{continuous_Hodges} have the same asymptotic properties of $\tilde\theta^{(i)}_{n}(c)$, $i=1,2$, as presented in Theorem \ref{multi-dimension-super-efficiency}.
\end{theorem}

\vskip3mm

\noindent{\bf 3. MLE version.}

\noindent
We next discuss the maximum likelihood estimation with the simplest i.i.d. case as an example. It is not difficult to extend the results to general situations. 
 
The log-likelihood function of $\theta$ from i.i.d. $X_1,\dots,X_n$ with a common density $f(x; \theta)$ is
 \[
 l(\theta)=l(\theta;X_1,\dots,X_n)=\prod_{i=1}^n \log f(X_i;\theta).
 \]
 It is well known that under certain regularity conditions, there exists a sequence of asymptotically efficient MLE $\hat\theta_n$, i.e., 
$\sqrt n(\hat \theta-\theta) \overset{d}{\rightarrow} N(0, {\frak F}^{-1}(\theta))$, where 
$${\frak F}(\theta)=-\E_\theta\left[{\partial^2\over\partial\theta\partial\theta'}\log f(X_1;\theta)\right]$$
is the Fisher's information matrix. For any $b\subset\{1,2,\dots,d\}$ and the corresponding rearrangement of $\theta=(\theta_b,\theta_{\bar b})$, ${\frak F}(\theta)$ can be rewritten as
 \[
 {\frak F}(\theta)=-\left(\begin{array}{cc}
\E_\theta\left[\displaystyle{\partial^2\log f(X_1;\theta)\over\partial\theta_b\partial\theta_b'}\right]
&\E_\theta\bigg[\displaystyle{\partial^2\log f(X_1;\theta)\over\partial\theta_b\partial\theta_{\bar b}'}\bigg]\\ 
 \E_\theta\left[\displaystyle{\partial^2\log f(X_1;\theta)\over\partial\theta_{\bar b}\partial\theta_b'}\right]
&\E_\theta\bigg[\displaystyle{\partial^2\log f(X_1;\theta)\over\partial\theta_{\bar b}\partial\theta_{\bar b}'}\bigg]
\end{array}\right)
=\left(\begin{array}{cc} {\frak F}_{bb'}(\theta)& {\frak F}_{b\bar b'}(\theta)\\
 {\frak F}_{\bar bb'}(\theta)& {\frak F}_{\bar b\bar b'}(\theta)\end{array}\right),\quad\hbox{say}.
 \]
For any constant vector $c=(c_1,\dots,c_d)$, Theorem \ref{multi-dimension-super-efficiency} yields the following immediate results.
  
\begin{theorem}\label{2.4}
 If $\{\hat\theta_n\}$ is an efficient sequence of maximum likelihood estimators, then the oracle Hodges' estimators $\tilde\theta(c)$ in Definition \ref{Hodges' estimators} have the following properties: 
For any $b\subset\{1,2,\dots,d\}$, $\theta_{\bar b}=c_{\bar b}$ implies 
$\lim_{n\rightarrow\infty}\limits\Pr(\tilde\theta_{n,\bar b}(c)=\theta_{\bar b})=1$ and
\[
\sqrt n(\tilde\theta(c)-\theta)\overset{d}{\rightarrow} N\left( 0 , 
\left(\begin{array}{cc} {\frak F}_{b(\theta) b'(\theta)}^{-1}(\theta)&0\\0&0\end{array}\right)\right).
 \]
 \end{theorem}
 
Clearly, $ {\frak F}_{bb'}(\theta)$ is the Fisher's information matrix for parameter $\theta_b$ depending on the unknown $\theta_{\bar b}$. The asymptotic variance of $\sqrt n(\hat\theta_{n,b}-\theta_{n,b})$ is 
 \begin{equation}\label{MLE_Fisher}
 ({\frak F}^{-1}(\theta))_{bb'}=({\frak F}_{bb'}(\theta)-{\frak F}_{b\bar b'}(\theta){\frak F}^{-1}_{\bar b\bar b'}(\theta){\frak F}_{\bar bb'}(\theta))^{-1}\geq {\frak F}^{-1}_{bb'}(\theta)
 \end{equation}
 with strict inequality if $b\neq \{1,2,\dots,d\}$.

Note that $ {\frak F}_{b(\theta) b'(\theta)}$ is the Fisher's information matrix of the marginal vector $\theta_{b(\theta)}$ knowing that $\theta_{\bar b(\theta)}=c_{\bar b(\theta)}$. Taking $b=b(\theta)$, Theorem \ref{2.4} shows the superefficiency of $\tilde\theta_n(c)$ over the MLE $\hat\theta_n$ at any $\theta$ in the parameter space such that $b(\theta)\neq \{1,2,\dots,d\}$ (i.e., $\theta_j=c_j$ for some $j\in\{1,2,\dots,d\}$). This result covers such parametric models as linear regression with normally distributed errors and generalized linear regression.
\vskip3mm

\newpage 

\noindent{\bf 4. LSE version.}

\noindent
For a linear regression model ${\bf Y}={\bf X}\beta+\varepsilon$ with $\E[\varepsilon]=0$ and Var$(\varepsilon)=\sigma^2I_n$, we can generate the oracle Hodges' estimator $\tilde\beta(c)$ from the least square estimate $\hat\beta_{LS}=({\bf X}'{\bf X})^{-1}{\bf X}'{\bf Y}$. Under certain regularity conditions (see e.g., van de Vaart, 2000, Example 2.28), $\hat\beta_{LS}$ is asymptotically distributed as $\sqrt n(\hat\beta_{LS}-\beta)\overset{d}{\rightarrow} N({\bf0}, \sigma^2\Sigma_X^{-1})$, where $\Sigma_X=\lim_{n\rightarrow\infty}\limits{n}^{-1}{\bf X}'{\bf X}$.
 By Theorem \ref{multi-dimension-super-efficiency} again, the following result is obvious.
 \begin{theorem}\label{LSE superefficiency}
Given any fixed $d$-vector $c$ and subset $b\subset\{1,2,\dots,d\}$,  $\beta_{\bar b}=c_{\bar b}$ implies that $\lim_{n\rightarrow\infty}\limits\Pr(\tilde\beta_{n,\bar b}(c)=\beta_{\bar b})=1$ and
 \begin{equation}\label{Sigma_b_beta}
 \sqrt n(\tilde\beta(c)-\beta)\overset{d}{\rightarrow}  N\left( 0 , \left(\begin{array}{cc}\sigma^2{\Sigma}_{b(\beta) b(\beta)}^{-1}&0\\0&0\end{array}\right)\right), \hbox{ where }{\Sigma}_{b(\beta) b(\beta)}=\lim_{n\rightarrow\infty}{1\over n}{\bf X}_{b(\beta)}'{\bf X}_{b(\beta)}. 
 \end{equation}
 \end{theorem}
This shows that $\sqrt n(\tilde\beta_{b(\beta)}(c)-\beta_{b(\beta)})$ has the same asymptotic distribution as the oracle estimator $\hat\beta_{b(\beta)}^o=({\bf X}_{b(\beta)}'{\bf X}_{b(\beta)})^{-1}{\bf X}_{b(\beta)}'(Y-{\bf X}_{\bar b(\beta)}c_{\bar b(\beta)})$ if the true value of $\beta$ is $(\beta_{b(\beta)}',c_{\bar b(\beta)}')'$.
\subsection{Two remarks}
We conclude this section with the following two remarks.
\begin{remark}\label{difference}
By  Theorems \ref{asymp_1}, the  classical version $\breve\theta_n(c)$  has an asymptotic distribution given by $r_{n}(\breve{\theta}_{n}(c)-\theta)\overset{d}{\rightarrow}ZI(\theta\not=c)$. That is, $\breve\theta_n(c)$ can only improve the asymptotic variance of $\hat\theta_n$ at $\theta=c$ in the parameter space $\Theta$, which is much more restrictive than the improvement achieved by the limit in \eqref{2.10} for $\tilde\theta_n(c)$ defined by \eqref{Hodges' estimate-multi-dimension}. To see this,  
note that the variance of $\check Z_{b(\theta)}$ is 
\[
\tilde V_{b(\theta)}=V_{{b}(\theta),{b}(\theta)}^{-1}(V_{{b}(\theta),{b}(\theta)}
\quad V_{b(\theta),\bar{b}(\theta)}) V^{-1} \left(\begin{array}{c}V_{b(\theta), b(\theta)}\\ 
V_{b(\theta),\bar{b}(\theta)}\end{array}\right) V_{b(\theta),b(\theta)}^{-1}
=V_{b(\theta),b(\theta)}^{-1}\leq \Delta_{b(\theta)},
\]
where $\Delta_{b(\theta)}$ is the asymptotic variance of $r_n(\hat\theta_{n,b(\theta)}-\theta_{b(\theta)})$ by \eqref{Delta_theta}, and the equality holds only when $\hat\theta_{n,b(\theta)}$ and $\hat\theta_{n,\bar{b}(\theta)}$ are asymptotically independent.
Therefore, $\tilde\theta_n(c)$ can improve the asymptotic variance of $\hat\theta_n$ at any $\theta$ with $\bar{b}(\theta)\neq\emptyset$. Note also that $r_n(\tilde\theta_n(c)-\theta)\rightarrow 0$ in probability at $\theta=c$. 
In terms of asymptotic variances, $\breve\theta(c)$ improves $\hat\theta_n$ and $\tilde\theta_n(c)$ further improves $\breve\theta(c)$. A further important feature of $\tilde\theta_n(c)$ is its ability to produce oracle model selection procedures due to its form of asymptotic variance, as we will show in the next subsection. In contrast, neither $\breve\theta_n(c)$ nor $\hat\theta_n$ has such a capacity. These together highlight  the significant differences between the new oracle Hodges' estimator $\tilde\theta_n(c)$ and the classical version $\breve\theta(c)$.
\end{remark}

\begin{remark}\label{uniform distribution}
Note that MLE and LSE and their versions of oracle Hodges' estimators are of root-$n$ consistency under relevant regularity conditions. But this is not necessary for Definition \ref{Hodges' estimators}. The general version of $\tilde\theta_n(c)$ and its continuous version do not require those regularity conditions and they are not necessarily of root-$n$ consistency. For example, let $Y_i=(Y_{i1},\dots,Y_{id})'$, $i=1,2,\dots,$ be independent with identical uniform distributions over $\prod_{k=1}^d[-\theta_k,\theta_k]$, $\theta_k>0$, so that the MLE of $\theta=(\theta_1,\dots,\theta_d)$ is $\hat\theta_{n}=(\hat\theta_{n1},\dots, \hat\theta_{nd})$ with $\hat\theta_{nk}=\max(|Y_{1k}|,\dots, |Y_{nk}|)$, $k=1,2,\dots,d$. Then $n(\hat\theta_{n}-\theta)\overset{d}{\rightarrow} Z=(Z_1,\dots,Z_d)$ with mutually independent components $Z_1,\dots,Z_d$ such that
\[
\Pr_\theta(Z_k\leq x_k)=\Bigg\{\begin{array}{ll}e^{x_k/\theta_k} &\hbox{if }x_k<0\\ 1 &\hbox{if }x_k\ge 0
\end{array},\quad k=1,2,\dots,d.
\]
In this case, the general version $\tilde\theta_n(c)$ in \eqref{Hodges' estimate-multi-dimension} and its continuous version $\tilde\theta_n(c,f)$ in \eqref{continuous_Hodges} are still valid even though the regularity conditions of the likelihood function are not satisfied, but they are not of root-$n$ consistency.
\end{remark}

%%%%%%%%%%%%%%%%%%%%%%%%%%%%%%%%%%%%%%%%%%%%%%%%%%%%%%%%%%%%%%

\section{Model selection function and oracle property of $\tilde\theta_n(c)$}\label{model selection}
Definition \ref{Hodges' estimators in multiple dimensions} and Theorem \ref{multi-dimension-super-efficiency} clearly indicate the following properties of $\tilde\theta_n(c)$:
\begin{itemize}
  \item[(1)] $\tilde\theta_n(c)$ is a sparse estimate in the sense that some components of $\theta$, say $\theta_j$,  may be estimated by component $c_j$ of $c$.
  \item[(2)] $\lim_{n\rightarrow\infty}\limits\Pr(\tilde\theta_{n,\bar b(\theta)}(c)=c_{\bar b(\theta)})=1$ for whatever true value of the parameter $\theta$.
  \item[(3)] For any sequence of estimators $\hat\theta_n$, it is possible to define a new sequence $\tilde\theta_n(c)$ such that its asymptotic covariance matrix (with the same convergence rate as $\hat\theta_n$) is 
\begin{itemize}
\item[$-$] equal to that of $\hat\theta_n$ if $\theta_j\neq c_j$ for all $j\in \{1,2,\dots,d\}$, i.e., 
$ b(\theta)=\{1,2,\dots,d\}$,
\item[$-$] positive definite and strictly less than that of $\hat\theta_n$ if 
$\emptyset\neq b(\theta)\neq \{1,2,\dots,d\}$; in this case $\tilde\theta_n(c)$ is asymptotically more efficient than $\hat\theta_n$ because  
\begin{equation*}\label{covHodges}
{\rm Cov}\left(\begin{array}{c}\check Z_{b(\theta)}\\
0\end{array}\right)=\left(\begin{array}{cc}V_{b(\theta)b'(\theta)}^{-1}&0\\0&0\end{array}\right)\leq V^{-1},
\end{equation*} 
\item[$-$] zero at $\theta=c$, i.e., $b=\emptyset$.
\end{itemize}
\end{itemize}

Taking the center parameter $c=0$, the estimator $\tilde\theta_n(0)$ obtained from any $\hat\theta_n$ provides a model selection procedure that removes any parameter $\theta_j$ estimated by $c_j=0$ from $\theta=(\theta_1,\dots,\theta_d)$. This allows any model that omits some or all $\theta_j$ to be selected. In contrast, the classical Hodges' estimator $\breve\theta(0)$ can only choose between two extreme models: the full model (corresponding to $\breve\theta(0)\neq0$) or the null model (corresponding to $\breve\theta(0)=0$), provided, with no loss of generality, that every element of the original $\hat\theta_n$ is nonzero. 

The model selection methods derived from $\breve\theta(0)$, $\tilde\theta(0)$, penalized maximum likelihood estimation (PMLE) and penalized least square estimation (PLSE) are discussed in more details below:
\begin{enumerate}
\item[(1)] While both Hodges' estimators $\breve\theta(0)$ and $\tilde\theta(0)$ can be applied as long as a good estimator $\hat\theta_n$ is available ($r_n$ consistent for some constant sequence $r_n\rightarrow\infty$), PMLE can only be applied when likelihood functions are available and PLSE is limited to regression models, both under certain regularity conditions (see, e.g., Fan and Li for a set of regularity conditions) to produce root-$n$ consistent estimators. This is demonstrated by the example discussed in Remark \ref{uniform distribution}, where both $\breve\theta(0)$ and $\tilde\theta(0)$ can be applied, but neither PMLE nor PLSE because the regularity conditions fail to hold.

The above comparison is limited to a fixed parameter space of dimension $p<n$. The case of $p>n$, with $p$ varying with $n$, is not  considered in this paper. It is noted that PMLE/PLSE to a large extent are motivated by the need to deal with $p>n$, and may have certain advantages in such a case. The performance of PMLE/PLSE, and whether there are better estimators than PMLE/PLSE for $p>n$, may be interesting subjects for further research, which is however beyond the intention and scope of the present paper.  

\item[(2)] If $\breve\theta_n(0)$ is used to select model by removing the parameters estimated by zero, then either all parameters are selected, or all excluded (cf. Remark \ref{difference}), so that the resulting model selection does not possess oracle property. This highlights a major difference between the new type of Hodges' estimators $\tilde\theta_n(c)$ and the classical version $\breve\theta_n(c)$. 

\item[(3)] Due to point (2) above, the model selection method based on $\breve\theta(0)$ is limited to two candidate models only: the null model and the full model.
%, which is hardly useful. 
In contrast, the other methods allow all submodels of the full model to be candidates with certain $\theta_i$ set to zero.

\item[(4)] The classical Hodges' estimators are a special type of preliminary-test estimators. It appears, however, that the finite sample behaviors of preliminary-test estimators have not been adequately addressed in the literature either, despite the importance of this problem in statistics. Some exceptions can be found in, e.g., the book by Judge and  Bock (1978) and the review paper of Giles and Giles (1993).
Although PMLE/PLSE and oracle Hodges' estimators may look like a preliminary-test estimation, there are fundamental differences:
\begin{itemize}
\item Existing preliminary-test procedures are essentially based on a single hypothesis that is either accepted or rejected as a whole.
\item The model selection procedures derived by PMLE/PLSE and oracle Hodges' estimation identify every single parameter in the multi-dimensional vector of parameters and judge if it is estimated by zero or some nonzero value. These are similar to estimation after multiple tests for a family of  hypotheses (multiple tests are also known as multiple comparisons, see., e.g., Hsu, 1996); it would be appropriate to call it preliminary-multiple-test estimation. In particular, the oracle Hodges' estimators provide an instance of such an estimation.
\end{itemize} 

\item[(5)] Properties (1)--(3) above ensure $\tilde\theta(0)$ model selection to possess the oracle property, provided the original estimator $\hat\theta_n$ is root-$n$ consistent and efficient -- which is the typical case where PMLE/PLSE can be employed to produce oracle model selection procedures and $\hat\theta_n$ is taken to be the efficient MLE/LSE. See, e.g., Fan and Li (2001) for a general discussion of the penalties producing oracle model selection procedures. As a result, properties (1)--(3) are more general than oracle properties, and both $\tilde\theta(0)$ and  PMLE/PLSE methods are sparse and have the same asymptotic behavior in the case of regular likelihood functions or regression models. 

\item[(6)] Note that, in particular, when $V$ is a diagonal matrix, the model selection driven by $\tilde\theta_n(0)$ is a hard-threshold model selection mentioned in Fan and Li (2001). If $V$ is not diagonal, then by Definition \ref{Hodges' estimators}, $\tilde{\theta}_j(0)=0$ for $j\in\bar{b}(0)$ and $\tilde\theta_j(0)$ is obtained from $\hat\theta_{nj}$ for $j\notin\bar{b}(0)$ with adjustments by all estimators $\hat\theta_{n1},\dots,\hat\theta_{nd}$, rather than solely from $\hat\theta_{nj}$. This also accounts for why the raw hard-thresholding model selection does not have the oracle property. 

\item[(7)] Leeb and P\"otscher (2005, 2008b) argued that SCAD and other sparse estimators are a return of Hodges' estimators by examining their finite sample performance in the neighborhoods of $\theta= 0 $. For any sparse estimator $\breve\theta$ satisfying $\lim_{n\rightarrow\infty}\limits\Pr_0(\breve\theta= 0 )=1$, where $\Pr_0$ is the probability computed at $\theta=0$, the maximum risk over a neighborhood of $\theta= 0$ tends to the maximum of the employed risk function as $n\rightarrow\infty$, and to infinity if the square loss is used. However, this section shows that SCAD or other oracle estimators are of the same asymptotic property with the new type of Hodges' estimator $\tilde\theta_n(0)$, rather than a simple return of the classical form $\breve\theta_n(0)$ because $\breve\theta_n(0)$ does not have the oracle property.

\item[(8)]  PMLE maximizes the penalized likelihood functions. It is however generally unknown whether a PMLE can be expressed as a function of MLE $\hat\theta_{ML}$. In fact, PMLE is  solved by direct maximization using numerical algorithms. In linear regression ${\bf Y}={\bf X}\theta+\varepsilon$, because  PLSE minimizes the penalized squares 
\[
({\bf Y}-{\bf X}\theta)' ({\bf Y}-{\bf X}\theta)+\sum_{i=1}^df_i(\theta_i,\lambda_n)={\bf Y}'(I-P){\bf Y}+(\theta-\hat\theta_{LS})'{\bf X}'{\bf X}(\theta-\hat\theta_{LS})+\sum_{i=1}^df_i(\theta_i,\lambda_n), 
\]
where $P$ is the projection matrix onto the column space of $\bf X$, it is clear that PLSE is a function of the LSE $\hat\theta_{LS}$. Its analytical form, however, is also generally unavailable and numerical methods are again needed to solve it. In contrast, if one takes $\hat\theta_n$ to be the MLE $\hat\theta_{ML}$ or LSE $\hat\theta_{LS}$ in these two cases, then both $\breve\theta(c)$ and $\tilde\theta(c)$ can be expressed explicitly by $\hat\theta_{ML}$ or $\hat\theta_{LS}$ as in Definitions \ref{Hodges' estimators_naive} and \ref{Hodges' estimators}. Due to these explicit forms of $\breve\theta(c)$ and $\tilde\theta(c)$, we can theoretically derive lower bounds for the performance of model selectors driven by $\breve\theta(c)$ or $\tilde\theta(c)$ in finite sample size in the next section.
\end{enumerate}

%%%%%%%%%%%%%%%%%%%%%%%%%%%%%%%%%%%%%%%%%%%%%%%%%%%%%%%%%%%%%%%%%%%%%%%%%%
\section{Performance of Hodges' estimators}\label{performance}
\setcounter{equation}{0}
Hodges' example gives a counterexample to the conjecture of R. A. Fisher that the MLE is asymptotically efficient with the smallest asymptotic variance (at rate $\sqrt n$).
To overcome the difficulty thus caused,  Le Cam (1953) and other researchers proved that
the improvement on MLE can only occur in a subset of the parameter space with zero Lebesgue measure. This argument, however, did not provide any reason to rule out the use of Hodges' estimators. An obvious question remains: if a Hodges' estimator can outperform the MLE at even one point without paying any price, why not use it as a preferred one? 

We now attempt to answer this question with the following two arguments:
\begin{itemize}
\item[(1)] Generally, the asymptotic bias and variance are not necessarily connected to their finite sample size versions, hence a small asymptotic variance does not imply a small variance of an estimator even if the sample size is very large. To link the asymptotics to the finite sample size, a further condition of uniform integrability is required.

\item[(2)] Even if the required uniform integrability is attached, the performance of Hodges' estimators $\tilde\theta_n(c)$ is still poor at the vicinity of $\theta=c$ due to a lack of uniformity in convergence over $\theta\in\Theta$. 
\end{itemize} 
For post-model-selection estimation in regression analysis, the lack of uniformity in convergence over parameters has been discussed by Yang (2005, 2007) and Leeb and P\"otscher (2008a, 2008b), among others.
In Subsection 4.1 below, we discuss this issue further with a more direct and elementary approach, which is made possible by the particular form of the oracle Hodges' estimators.

\subsection{Uniform integrability and asymptotics}

First we recall some misconception regarding  asymptotic mean-squared error (MSE) and asymptotic efficiency  of an estimator. Given a normalized sequence $r_n(\hat\theta_n-\theta)\overset{d}\rightarrow Z$, the asymptotic mean and variance (hence MSE) are only the mean and variance of $Z$, but not the limit of $r_n\E[\hat\theta_n-\theta]$ and $r_n^2{\rm Var}(\hat\theta_n-\theta)$ in general. 

It is an easy exercise to construct examples in which a sequence of unbiased (biased) estimates might be asymptotically biased (unbiased), and even in the class of unbiased estimators, a sequence of estimates with smaller variances for every sample size $n$ might have larger asymptotic variance, and vice versa. Consequently, the concepts of asymptotically unbiased estimate and asymptotic variance could be highly misleading when they are considered as analogies to unbiased estimate and variance in finite samples. Unfortunately, this lesson appears to have been overlooked by many statisticians for a long time, especially when one proves the oracle property of an estimator.
 
Let $\hat\theta_n$ be a sequence of estimates such that $r_n(\hat\theta_n-\theta)\overset{d}{\rightarrow} Z$, where $r_n\rightarrow\infty$ as $n\rightarrow\infty$ and $Z$ is a random variable with mean zero and finite variance $\sigma^2$.
In this context, for $Z_n=\hat\theta_n-\theta$ and $Y_n=r_n(\hat\theta_n-\theta)$, we have the following two facts by Theorem 2.20 of van der Vaart (1998):
 \begin{itemize}
 \item[(i)]  $\E[Z_n]\rightarrow0$ ($\E[Y_n]\rightarrow0$) if and only if $\{Z_n:n\geq1\}$ ($\{Y_n:n\geq1\}$) is uniformly integrable.
 \item[(ii)] $\E[Z_n^2]\rightarrow0$ ($\E[Y_n^2]\rightarrow \sigma^2$) if and only if $\{Z_n^2:n\geq1\}$ ($\{Y_n^2:n\geq1\}$) is uniformly integrable.
 \end{itemize}
 
 Obviously, conditions like these have been completely neglected in the literature seeking estimators possessing the oracle property.
 %%%%%%%%%%%%%%%%%%%%%%%%%%%%%%%%%%%%%%%%%%%%%%%%%%%%%%%%
\subsection{Performance of Hodges' estimators}
To ensure the finite sample size quantities to approach their asymptotic versions as the sample size increases, additional conditions are required.  In the remainder of this section, we assume that the required uniform integrability described above is satisfied, so that the bias and variance in finite sample size approach their asymptotic versions when the sample size is large. Even in such a case, however, both numerical and theoretical analyses below show that the overall performance of Hodges' estimators in finite sample does not match its limit as analyzed in this subsection below.

Due to the difficulty to obtain the exact MSE of Hodges' estimator in closed form, earlier arguments against the use of Hodges' estimators were largely based on numerical results. For example, a result for Hodges' estimate of the mean $\theta$ in the normal distribution $N(\theta, 1)$ is well-known in the literature, see e.g., van der Vaart (1998) and Lehmann and Casella (1998), which is recalled here.  
Let $\hat\theta_n=\bar{X}\sim N(\theta, {1/n})$. Taking $c=0$ and $a_n=n^{-1/4}$ yields the original version of Hodges' estimator $\tilde\theta_n(0) = \bar{X}I(|\bar X|>n^{-1/4})$. Note that the sequences $\sqrt n(\hat\theta_n-\theta)$, $\sqrt n(\tilde\theta_n(0)-\theta)$,  $ n(\hat\theta_n-\theta)^2$ and  $n(\tilde\theta_n(0)-\theta)^2$ are all uniformly integrable so that the asymptotic means and variances are equal to the respective limits of the means and variances in finite sample size.
While the MSE (scaled by $n$) of $\bar{X}$ is constant $1$ for all $\theta$'s, that of the Hodges' estimator $\tilde\theta_n(0)$ can only be numerically computed. Figure 1 shows the curves of the MSE of $\tilde\theta_n(0)$ in $\theta$ for sample sizes $n=5$, 50 and 500, which behave poorly in the vicinity of zero, particularly at large sample size ($n=500$). This illustrates a much worse performance of the Hodges' estimator $\tilde\theta_n(0)$ than the MLE $\hat\theta_n$ near the center $c=0$.

\begin{figure}[h]
\caption{The MSE (scaled by $n$) of Hodges' estimator\label{Hodges}}
\begin{center} 
{\includegraphics[width=3.5in, height=2.8in]{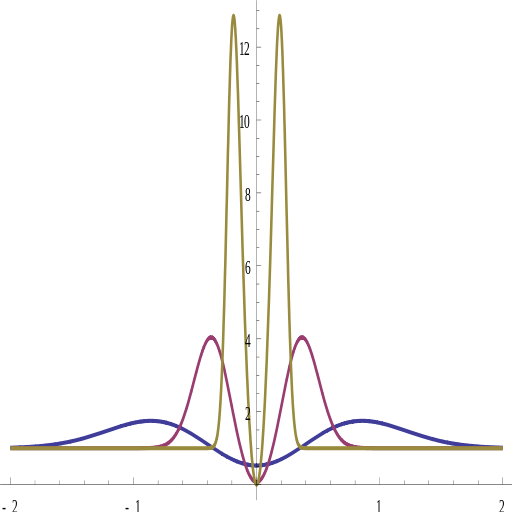}}\\
%{\includegraphics[natheight=7.5cm, natwidth=11.5cm, height=5cm, width=8.5cm]{Hodges_estimator_risk_function.png}}\\
Legend: Blue curve corresponds to $n = 5$, purple to $n = 50$, and olive to $n = 500$.
\end{center}
\end{figure}

Although finding the exact closed form of the MSE (or other performance measures) of a Hodges' estimator is difficult, even in the simple case discussed just now, it turns out that some useful lower bounds of regular losses of Hodges' estimators $\tilde\theta_n(c)$ can be obtained to see the rationale behind their poor performance in the vicinity of the center point $c$ as well as of the subsets of the parameter space with $\theta_j=c_j$ for some $j\in\{1,2,\dots,d\}$. This seems to have been overlooked by the community. 

The results presented below are from theoretical analyses on a general ground and distribution-free -- they are valid regardless of the underlying distributions of the population. 

\vskip4mm
\noindent {\bf Performance of classical Hodges' estimators}
\vskip2mm

\noindent
We first present the results for classical Hodges' estimators, which are in fact almost sure results.
   
\begin{theorem}\label{prop_Hodges}
    Under the conditions of Theorem \ref{asymp_1}, given any $k>0$, there exists a (deterministic) $N>0$ such that for all $n>N$,
     \begin{equation}\label{MSE_Hodges}
     r_n\|\breve\theta_n(c)-\theta\|\geq k\quad\hbox{for all }\theta \hbox{ satisfying }k\leq r_n\|\theta-c\|\leq a_nr_n-k.%\in\bigg[c-a_n+{k\over r_n},c-{k\over r_n}\bigg]\bigcup\bigg[c+{k\over r_n},c+a_n-{k\over r_n}\bigg]
     \end{equation} 
Furthermore, for $\theta_n$ such that $\|\theta_n-c\|= a_n/2$,
\begin{equation}\label{maximumMSE_Hodges}
     r_n\|\breve\theta_n(c)-\theta_n\|\geq{1\over 2}r_na_n\rightarrow\infty\quad
\hbox{as }n\rightarrow\infty. 
     \end{equation}
     \end{theorem}
     \begin{m_proof}
Since $a_n\rightarrow0$ and $a_nr_n\rightarrow\infty$ as $n\rightarrow\infty$ under the conditions of Theorem \ref{asymp_1}, the set $\{\theta: k\leq r_n\|\theta-c\|\leq a_nr_n-k \}$ is a nonempty ring when $n$ is sufficiently large such that $a_nr_n>2k$. For any $\breve\theta_n(c)$ defined in \eqref{Hodges_one_D}, for any values the sample may take, as long as 
    $\theta$ satisfies $ k\leq r_n\|\theta-c\|\leq a_nr_n-k$, i.e.,  $ k/r_n\leq \|\theta-c\|\leq a_n-k/r_n$, it is clear that 
    \begin{align*}
    \|\breve\theta_n(c)-\theta\|=&\|c-\theta\|I(\|\hat\theta_n-c\|\leq a_n)+\|\hat\theta_n-\theta\|I(\|\hat\theta_n-c\|> a_n)\\
   \geq& \|c-\theta\|I(\|\hat\theta_n-c\|\leq a_n)+(\|\hat\theta_n-c\|-\|c-\theta\|) I(\|\hat\theta_n-c\|> a_n)\\
   \geq &k/r_nI(\|\hat\theta_n-c\|\leq a_n)+(a_n-(a_n-k/r_n) I(\|\hat\theta_n-c\|> a_n)\\
   = &{k/r_n}.
    \end{align*}
%$\theta\in[c-a_n+{k/r_n},c-{k/r_n}]\cup[c+{k/r_n},c+a_n-{k/r_n}]$.
     This proves the first assertion in \eqref{MSE_Hodges}.
    The second assertion \eqref{maximumMSE_Hodges} is obvious. 
     \end{m_proof}
\vskip2mm

The results of Theorem \ref{prop_Hodges} have the following easy implications:
\begin{enumerate}
\item[(1)] Formula \eqref{maximumMSE_Hodges} proves that even if the MSE of $\breve\theta_n(c)$ (scaled by $r_n$) converges to that of the asymptotic distribution, the convergence is not generally uniform because
\[
r_n^2\max_{\theta\in\Theta}\E_\theta\big[\|\breve\theta_n(c)-\theta\|^2\big]
\geq r_n^2\E_{\theta_n}\big[\|\breve\theta_n(c)-\theta_n\|^2\big]
\geq{1\over 4}r_n^2a_n^2\rightarrow\infty.
\]
    
 \item[(2)] The same arguments also hold for more general loss functions $l(\hat\theta_n; \theta)=l(\|\hat\theta_n-\theta\|)$ with some nondecreasing function $l(u)$ in $u>0$ satisfying $l(0)=0$, so that the risk function scaled by $1/l(1/r_n)$ is $R_n(\hat\theta_n; \theta)=\E_\theta[{l(\|\hat\theta_n-\theta\|)/l(1/r_n)}]$. Let $p=\min\{i\geq 1:l^{(i)}(0+)\neq0\}$. %The results are stated in the next theorem.     
%\begin{theorem}
Then under the conditions of Theorem \ref{asymp_1}, for any $k>0$, there exists $N>0$ such that for all $n>N$,
\[
{l(\breve\theta_n(c);\theta)\over l(1/r_n)}\geq {l(k/r_n)\over l(1/r_n)}\rightarrow k^p \hbox{ as } n\rightarrow\infty \
\] 
for all $\theta$ satisfying $k\leq r_n\|\theta-c\|\leq a_nr_n-k$. Furthermore, for all $\theta_n$ such that $\|\theta_n-c\|=a_n/2$,
\[
{l(\breve\theta_n(c);\theta_n)\over l(1/r_n)}\geq{l(a_n/2)\over l(1/r_n)}\rightarrow\infty \hbox{ as } n\rightarrow\infty.
\]
The last formula also implies 
$$\max_{\theta\in\Theta}R_n(\breve{\theta}_n(c);\theta)\geq R_n(\breve{\theta}_n(c);\theta_n)\geq {l(a_n/2)\over l(1/r_n)}\rightarrow\infty \hbox{ as } n\rightarrow\infty.$$ 

\item[(3)] Another way is to analyze a loss function $L(\hat\theta_n; \theta)=l(r_n\|\hat\theta_n-\theta\|)$ with a nondecreasing function $l(u)$ in $u>0$ (cf. Leeb and P\"{o}tscher (2008b)), which corresponds to a sequence of loss functions $l_n(u)=l(r_nu)$ ($r_n=\sqrt{n}$ in their paper), so that 
$R_n (\hat\theta_n;\theta)=\E_\theta[L(\hat\theta_n; \theta)]$. 

Then similarly under the conditions of Theorem \ref{asymp_1}, given any $k>0$, there exists $N>0$ such that
$L(\breve\theta_n(c);\theta)\geq l(k)$ for all $n>N$ and $\theta$ satisfying $ k\leq r_n\|\theta-c\|\leq a_nr_n-k$.
Moreover, for all $\theta_n$ such that $\|\theta_n-c\|=a_n/2$,
$L(\breve\theta_n(c);\theta_n)\geq l({r_n a_n/2})\rightarrow l(\infty)$ as $n\rightarrow\infty$. 
In particular, the last property indicates that
\begin{equation}\label{Leeb and Potscher 2008b}
\max_{\theta\in\Theta}R_n(\breve\theta_n(c);\theta)\geq \E_{\theta_n}[L (\breve\theta_n(c);\theta_n)]\rightarrow l(\infty).
\end{equation}
This shows that even if $\lim_{n\rightarrow\infty}\limits R_n(\breve\theta_n(c);\theta)=R(Z;\theta)=\E_\theta[l(Z;\theta)]$ for every $\theta$ pointwise, the maximum risk over $\theta$ may increasingly tend to $l(\infty)$. An example is $l(u)=I_{(z,\infty)}(u)$ for any fixed continuity point $z\in\mathbb{R}^+$ of the distribution of $Z$. The risk function of an estimator $\hat\theta_n$ with this $l(u)$ is
$R_n(\hat\theta_n;\theta)=\Pr_\theta(r_n\|\hat\theta_n-\theta\|> z)$. Thus 
$$\lim_{n\rightarrow\infty}\limits R_n(\breve\theta_n(c);\theta)=\Pr_\theta(\|Z\|>z)I(\theta\neq c).$$
In contrast, $\Pr_\theta(r_n\|\breve\theta_n(c)-\theta\|>z)=1$ for all $\theta$ with $ k\leq r_n\|\theta-c\|\leq a_nr_n-k$ and $z\in[0,k]$, regardless how large is $n$. Moreover, for any $z>0$, if $\theta_n=\pm {a_n/2}$ and $n$ is sufficiently large such that $a_nr_n>x$, then $\Pr_{\theta_n}(r_n\|\breve\theta_n(c)-\theta_n\|>z)=1$. Consequently, 
\[
 \max_{\theta\in\Theta}R_n(\hat\theta_n;\theta)=\max_{\theta\in\Theta}\Pr_\theta(r_n\|\breve\theta_n(c)-\theta\|>z)=1.
 \]
 \end{enumerate}
\vskip4mm

\noindent {\bf Performance of the oracle Hodges' estimators}

\noindent
Now we turn to  analyze the performance of the oracle Hodges' estimators $\tilde\theta_n(c)$ defined in Section \ref{Hodges' estimators in multiple dimensions}.  To simplify the exposition, we assume without loss of generality that the parameter space is $\Theta=\mathbb{R}^d$ and define the following subsets of $\Theta$: 
  \begin{align*}
   &\Theta_{n1}=\{\theta: \min_{1\leq j\leq d}|\theta_{nj}-c_j|>a_{nj}\},&\bar\Theta^k_{n1}=\{\theta: d(\theta, \Theta_{n1})\leq {k/r_n}\},\\
   &\Theta_{n2}=\{\theta:\min_{1\leq j\leq d}|\theta_{nj}-c_j|=0\},  
&\bar\Theta^k_{n2}=\{\theta: d(\theta, \Theta_{n2})\leq {k/r_n}\},\\
%   &\Theta_{n3}=\Theta-\Theta_{n1}(c)\cup\Theta_{n2}(c) ,
&\bar\Theta^k_{n3}=\Theta-\bar\Theta^k_{n1}(c)\cup\bar\Theta^k_{n2}(c),\ \ 
   \end{align*}
which are all dependent on the center point $c$.
%Let  $\|\theta\|=\sqrt{\sum_{j=1}^d\theta_{j}^2}$ denote the Euclidean norm of $\theta\in\Theta=\mathbb{R}^d$. 
Under the uniform integrability of $\hat\theta_n$, Theorem \ref{multi-dimension-super-efficiency} indicates that
    \[
   \lim_{n\rightarrow\infty}r^2_n\E\big[\|\tilde\theta_n(c)-\theta\|^2\big]
\leq\lim_{n\rightarrow\infty} r^2_n\E\big[\|\hat\theta_n-\theta\|^2\big]
     \]
with the strict inequality at certain values of $\theta$. For fixed sample size $n$, however, we have the following contrary results, which are  extensions of Theorem \ref{prop_Hodges}.

\begin{theorem}\label{multi-bounds}
Under the conditions of Theorem \ref{multi-dimension-super-efficiency}, for any $k>0$, there exists $N>0$ such that
\begin{equation} \label{MSE_Hodge_multi}
r_n\|\tilde\theta_n(c)-\theta\|\geq k\hbox{ for all }\theta\in\bar\Theta^k_{n3}\quad\forall n>N,
\end{equation}
and for any $\theta^{(n)}\in\bar\Theta_{n3}^k$ satisfying $\theta^{(n)}_j=c_j\pm a_{nj}/2$,
\begin{equation}\label{maxMSE_Hodge_multi}
r_n\|\tilde\theta_n(c)-\theta^{(n)}\|
\geq {1\over 2}r_na_{nj}\rightarrow\infty\quad\hbox{as }n\rightarrow\infty.
\end{equation}
\end{theorem}

\begin{m_proof}
 Again, we assume that $N$ is sufficiently large such that $r_n\max_{1\leq j\leq d}a_{nj}>k$ for all $n>N$. It can be easily shown that $d(\Theta_{n1},\Theta_{n_2})=\min_{1\leq j\leq d}a_{nj}$. Because for any $\theta_1\in \bar\Theta^k_{n1}\hbox{ and }\theta_2\in\bar\Theta^k_{n2}$, 
    \[
    d(\Theta_{n1},\Theta_{n_2})\leq d(\Theta_{n1},\theta_1)+d(\theta_1,\theta_2)+d(\theta_2,\Theta_{n_2})\leq {2k\over r_n}+d(\theta_1,\theta_2),
    \]
the condition 
$\lim_{n\rightarrow\infty}\limits r_n\min_{1\leq j\leq d} a_{nj}=\infty$ implies
  \[
  d(\bar\Theta^k_{n1},\bar\Theta^k_{n2})=\min_{\theta_1\in \bar\Theta^k_{n1},\theta_2\bar\Theta^k_{n2}}d(\theta_1,\theta_2)\geq \min_{1\leq j\leq d}a_{nj}-{2k\over r_n}>0
  \]
for sufficiently large $n$.  Consequently,  $\bar\Theta_{n1}^k\cap\bar\Theta_{n_2}^k=\emptyset$ and $\bar\Theta^k_{n3}\neq\emptyset$.
  
By Definition \ref{Hodges' estimators} of $\tilde\theta_n(c)$,  $\hat\theta_n\in\Theta_{n1}$ implies $\tilde\theta_n(c)=\hat\theta_n\in\Theta_{n1}$ and $\hat\theta_n\in\Theta-\Theta_{n1}$ implies $\bar b_n(c)\neq \emptyset$, so that $\tilde\theta_n(c)\in\Theta_{n2}$ because $\tilde\theta_{n,\bar b_n(c)}(c)-c_{\bar b_n(c)}=0$. That is, $\tilde\theta_n(c)$ takes values only in $\Theta_{n1}\cup\Theta_{n2}$.
For any $\theta\in\bar\Theta^k_{n3}$, it is clear that $d(\theta,\Theta^k_{n1}\cup\Theta^k_{n2})\geq{k/r_n}$. 
As a result, $\theta\in\bar\Theta^k_{n3}$ implies 
\[
\|\tilde\theta_n(c)-\theta\|\geq d(\theta, \Theta^k_{n1}\cup\Theta^k_{n2})\geq{k\over r_n}.
\] 
Thus the assertion in \eqref{MSE_Hodge_multi} follows.
The second assertion in \eqref{maxMSE_Hodge_multi} is easy to check.
\end{m_proof}

\bigskip
The following extensions of Theorem \ref{MSE_Hodge_multi} are minor modifications of the points presented earlier for classical Hodges' estimators.

\begin{enumerate}
\item[(1)] Formulas \eqref{MSE_Hodge_multi} and \eqref{maxMSE_Hodge_multi} prove that, under the conditions in Theorem \ref{multi-dimension-super-efficiency}, for any $k>0$, there exists $N>0$ such that
\[
r_n^2\E_\theta\big[\|\tilde\theta_n(c)-\theta\|^2\big]\geq k^2\quad\hbox{ for any }\theta\in\bar\Theta^k_{n3}
\hbox{ if } n>N
\]
and
\[
\max_{\theta\in\Theta}r_n^2\E_\theta\big[\|\tilde\theta_n(c)-\theta\|^2\big]
\ge {1\over 4}r_n^2a_{nj}^2\rightarrow\infty\quad\hbox{as }n\rightarrow\infty.
\]
\item[(2)] If we use the loss function $L(\hat\theta_n;\theta)=l(\|\hat\theta_n-\theta\|)$ with a nondecreasing function $l(u)$ in $u>0$, so that the risk function scaled by $1/l(1/r_n)$ is $R_n(\hat\theta_n; \theta)=\E_\theta[{l(\|\hat\theta_n-\theta\|)/l(1/r_n)}]$, then for sufficiently large $n$, 
\[
{l(\tilde\theta_n(c);\theta)\over l(1/r_n)}\geq {l(k/r_n)\over l(1/r_n)}\rightarrow k^p \quad\hbox{as } n\rightarrow\infty 
\quad\hbox{for all }\theta\in\bar\Theta^k_{n3};
\]
and for $\theta^{(n)}$ with $\theta^{(n)}_j=c\pm a_{nj}/2$ for some $j$,
\[
{l(\tilde\theta_n(c);\theta^{(n)})\over l(1/r_n)}\geq{l(a_{nj}/2)\over l(1/r_n)}\rightarrow\infty \quad
\hbox{as } n\rightarrow\infty.
\]
The last formula also implies that
$$\max_{\theta\in\Theta}R_n(\tilde{\theta}_n(c);\theta)\geq R_n(\tilde{\theta}_n(c);\theta^{(n)})\geq {l(a_{nj}/2)\over l(1/r_n)}\rightarrow\infty \quad \hbox{as } n\rightarrow\infty.$$

\item[(3)] If we analyze a loss function $L(\hat\theta_n; \theta)=l(r_n\|\hat\theta_n-\theta\|)$ with a nondecreasing function $l(u)$ in $u>0$, as in Leeb and P\"{o}tscher (2008b), so that $R_n (\hat\theta_n;\theta)=\E_\theta[L(\hat\theta_n; \theta)]$, %The corresponding result is presented in the theorem below.
%\begin{theorem}\label{cor_Hodges}
then under the conditions of Theorem \ref{multi-dimension-super-efficiency}, for any given $k>0$, there exists $N>0$ such that for all $n>N$, $L(\tilde\theta_n(c);\theta)\geq l(k)\ \hbox{for all }\theta\in\bar\Theta_{n3}^k$; 
and $L(\tilde\theta_n(c);\theta^{(n)})\geq l({r_n a_{nj}/2})\rightarrow l(\infty)$
as $n\rightarrow\infty$ for $\theta^{(n)}$ with $\theta^{(n)}_j=c\pm a_{nj}/2$ for some $j$.
The last property implies 
\begin{equation}\label{Leeb and Potscher 2008b}
\max_{\theta\in\Theta}R_n(\tilde\theta_n(c);\theta)\geq \E_{\theta_n}[L (\tilde\theta_n(c);\theta_n)]\rightarrow l(\infty).
\end{equation}
Because of this fact, even if $\lim_{n\rightarrow\infty}\limits R_n(\tilde\theta_n(c);\theta)=R(Z;\theta)=\E_\theta[l(Z;\theta)]$ for every $\theta$ pointwise, the maximum risk over $\theta$ may increasingly tend to $l(\infty)$. 
\end{enumerate}
 
 \begin{remark}
We conclude this section by the following two points that highlight the difference between our work and those of others, such as Leeb and  P\"{o}tscher (2008b). 
\begin{enumerate}
\item[(1)] Leeb and P\"{o}tscher (2008b) proved a result similar to equation \eqref{Leeb and Potscher 2008b}, which is more general with an arbitrary sparse estimator (say, $\tilde\theta_n$), but restricted to the regular case of $r_n=\sqrt{n}$, $c=0$ and normally distributed $Z$. It also requires the condition that $P_{n,k/\sqrt{n}}$ is contiguous with respect to $P_{n,0}$, where $P_{n,\theta}$ is the distribution of $\tilde\theta_n$, or a stronger condition that $P_{n,\theta}$ is locally asymptotically normal.  In comparison, we obtained stronger results expressed in \eqref{MSE_Hodge_multi} and \eqref{maxMSE_Hodge_multi} for classical and oracle Hodges' estimators,  which hold almost surely without such conditions as $\sqrt{n}$-consistency and contiguity. %This remark also applies for multi-dimensional parameters.
\item[(2)] Leeb and  P\"{o}tscher revealed the erratic behavior of a sparse estimator only in the vicinity of $0$, whereas we here showed that the erratic behavior of the oracle model selection procedure derived by a Hodges' estimator occurs not only in the vicinity of the point $c$ but also in the vicinity of every subset of $\Theta$ with some component $\theta_j=c_j$.
\end{enumerate}
\end{remark}

\section{Further discussions}\label{concluding}
 
To sum up, we have in this paper demonstrated that:
\begin{itemize}
\item[(1)] the oracle model selection procedures are not a simple return of the classical Hodges' estimators but more like \textit{oracle Hodges' estimators} in asymptotic sense;
\item[(2)]  properly constructed oracle Hodges' estimators can easily generate oracle model selection procedures that satisfy the requirements of continuity or smoothing;
\item[(3)] under the MSE criterion, the oracle Hodges' estimator $\tilde\theta_n(c)$ does not outperform its origin $\hat\theta_n$ in finite sample size, despite having a smaller asymptotic variance than $\hat\theta_n$; and 
\item[(4)] Hodges' estimators possessing the oracle property can perform much worse than their origins under the minimax criterion. 
\end{itemize}

Similar results to (3) and (4) can be found in Yang (2005, 2007) and the works by Leeb and his collaborators as mentioned above. The difference here is that the particular form of the oracle Hodges' estimators makes the proof quite obvious.

Points (3) and (4) above also provide an answer to the question why Hodges' estimators are not preferable to use even if they can improve the asymptotic efficiency at no cost. The key reasons behind this answer, as discussed earlier, are the disconnection between the performances of finite sample statistics and their asymptotics in certain situations (principally due to the lack of uniformity in integrability of the statistics), and  the universal lack of uniformity in the convergence in the situations where even the convergence is guaranteed.

In addition, a few points worth for further discussions are listed below.

\bigskip

\noindent \textbf{1. Uniformity in integrability and the convergence}.

\medskip\noindent 
By ignoring the uniformity in integrability and convergence, the widely adopted concept of asymptotically efficient estimation tends to place too much weight on the asymptotic distribution, leading to inappropriate use of asymptotic variance to measure the goodness of estimators. This creates the concept of superefficiency and supports Hodges estimators as superior, despite their poor performance in finite sample size. Because Hodges' estimators can be easily applied to generate model selection procedures possessing the oracle property, as demonstrated in point (2) above, our results also suggest that the oracle property itself is not a convincing reason to recommend the application of such oracle model selection procedures in theoretical studies or real-life applications.
  
This phenomenon was in fact recognized by statisticians decades ago. See for example H\'ajek (1971, p. 153)'s warning quoted below (which can also be found in Hu\v{s}kov\'a, Beran and Dupac, 1998, pp. 613-614):
\begin{quote} 
``Having obtained an asymptotic result we are not usually able to tell how far it applies to particular cases with finite $n$. $\dots$. Consequently, in applications we are guided by two epistemologically very different knowledge: (i) we have limit theorems giving some hope, but not assurance, of practical sample sizes; (ii) we work with some numerical experience, which we extend to cases that seem to us to be similar. $\dots$. 
Especially misinformative can be those limit results that are not uniform. Then the limit may exhibit some features that are not even approximately true for any finite $n$.
$\dots$. Superefficient estimates produced by L.J. Hodges (see Le Cam (1953)) have their amazing properties only in the limit. For any finite $n$ they behave quite poorly for some parameter values. These values, however, depend on $n$ and disappear in the limit.''
\end{quote}
Therefore, in the cases relying on asymptotic distribution of a statistic, uniform integrability and convergence of the statistic are important and desired properties to ensure the proximity between the finite-sample version and its asymptotics, especially when the statistic does not have an analytical form and can only be obtained by numerical computations.

\bigskip
\noindent \textbf{2. The power of simulations}.

\medskip\noindent 
Undoubtedly, with the development of contemporary computing facilities, Monte Carlo simulations provide more and more accessible and powerful tools in exploring  properties of statistical inferences. This powerful tool, however, appears being applied excessively to justify certain theoretical properties that are difficult to prove mathematically. It is even a common practice that, instead of providing theoretical results on  finite sample behaviors, many statisticians justify their statistical methods by finding and theoretically proving  asymptotic properties of their statistical methods and then supporting their finite sample  properties by Monte Carlo simulations. 

Therefore, another point worth to mention is the need to avoid the pitfall of relying too much on simulations to judge the merits of a property or a method, particularly in the era of fast advancing computer technology and capacity, as demonstrated by Leeb and P\"otscher (2008b) on the SCAD estimator in Example 4.1 of Fan and Li (2001). This pitfall is generally recognized in the statistical community, but appears often overlooked or ignored in a large volume of publications in the pursuit of finding new and exciting ideas and methods. 
Logically, simulation results can be convincing to counter a claim or conjecture made on a general ground, but not to support such a claim or conjecture, because it is impossible to exhaust all possible scenarios by simulations.  An assertion on property or goodness of a statistical procedure should be theoretically examined if it is to be claimed on a general ground, rather than rely on limited numerical simulations only.
\bigskip

\noindent \textbf{3. Generality of the center point}.

\medskip\noindent
The generality of the center point $c=(c_1,\dots,c_d)$ implies that, mathematically, one can make model selection not only on pre-specified $\theta_i$ but also on more general $\theta_i-c_i$, $1,\dots,d$. Take the oracle Hodges' estimators as an example, for any given constant $c=(c_1,\dots,c_d)$, one can make a model selection mathematically with the collection of candidate models defined by setting $\theta_i=0$ for one or more $i\in\{1,\dots,d\}$. Similarly for PMLE/PLSE, the common practice is to use penalty $f_i(|\theta_i|, \lambda_n)$ to select the ``best'' with certain $\theta_i=0$ from the candidates. By the same method and logic, for any fixed $c=(c_1,\dots,c_d)$, it is equally feasible to use penalty $f_i(|\theta_i-c_i|,\lambda_n)$ to select the ``best'' model with $\theta_i=c_i$ for certain $i$. This raises an obvious question in practice: why does everyone only take $c=0$ for model selection, not any other $c$?  As an example, if $m$ analysts under a common supervisor independently analyze a same data set with parameter $\theta$ of interest using oracle Hodges' estimation, PMLE or PLSE with different points $c^{(1)},\dots,c^{(m)}$, including zero or not, then they obtain $m$ different estimates $\hat\theta^{(j)}$, $j=1,\dots,m$. Clearly, all these estimators are of similar statistical properties and rationales, but there is no statistically sensible guidance for the supervisor to decide which one should be used. This question appears difficult to answer both logically and philosophically, and it raises a sobering question on the validity of the commonly adopted model selection procedures.
\bigskip

\noindent\textbf{4. Model selection for predictions.}

\medskip\noindent
In real applications of regression analysis, as well as in popular textbooks, variable selection is also discussed with a purpose of prediction, or equivalently, estimation of the regression function at certain points. In this paper, our efforts are focussed on the performance of the parameter estimators in model selection. This does not, however, point to a lack of generality of our results, due to the following two reasons:
\begin{itemize}
\item[(1)] In the case of linear regression with predictors represented as a $p$-vector ${\bf x}$ (whether the intercept term $1$ is included or not is irrelevant), taking $\theta$ as the vector of regression coefficients. Suppose that the regression function is to be estimated at $d$ point $\tilde{\bf x}_1,\dots,\tilde{\bf x}_d$ under the performance measure $A=\E\|\tilde{\bf X}\hat\theta-\tilde{\bf X}\theta\|^2=\sum_{i=1}^d\E(\tilde{\bf x}_i'\hat\theta-\tilde{\bf x}_i'\theta)^2$, where $\tilde{\bf X}=(\tilde{\bf x}_1,\dots,\tilde{\bf x}_d)'$ collects the data of the predictors at the $d$ prediction  points. Then the fundamental properties of positive definite matrices state that the order of different estimates under measure $A$ can be derived from the order of the covariance matrices of the corresponding estimates, given the presence of the latter, in the sense that for two matrices $C$ and $D$, $C>D$ if and only if ${\bf x}'C{\bf x}>{\bf x}'D{\bf x}$ for any dimension-compatible vector ${\bf x}$. As a result, the performance of the parameter estimators in our discussions is equivalent to that of the regression function estimator.
\item[(2)] In case the regression function is possibly nonlinear in the parameters $\theta$, including generalized linear models and general parametric nonlinear regression models, the similar arguments work asymptotically with the help of linear approximation of  statistics (i.e., the commonly known delta-method).

\end{itemize}

\noindent \textbf{5. Open questions}.

\medskip\noindent
Back to the controversy on the merits of the oracle property and the efforts to find oracle model selection procedures, the poor performance of Hodges' estimators and their oracle property appear to support the criticisms of the oracle estimators and cast serious doubts on the usefulness of oracle procedures. This further casts doubts in the validity and usefulness of prevailing model selection methods. We believe that the following open and challenging questions need to be convincingly answered before a consensus can be reached one way or the other: 
\begin{itemize}
\item[(1)] Without counting on the oracle property, are there any theoretical properties that ensure good performance of the popular model selection procedures, such as LASSO and SCAD? For example, is there any reasonable loss function beyond the regular ones discussed in this paper such that some of the popular model selection procedures can perform well? 
\item[(2)] Our arguments are limited to the case where the MLE/LSE exists so that $n>d$ is implicitly assumed. The cases of $d>n$, $d>>n$, or more generally, the design matrices of reduced-ranks, have so far been discussed under the assumption of parsimony with few exception, which are significantly different from the classical model. We have not yet come up with a clear idea on how the oracle model selection performs under this significantly different setting in  finite sample size. Further efforts are needed to examine various important topics regarding $d>n$ or $d>>n$, including model selection, which will be in the agenda of our future research works.
\end{itemize}
These open questions call for further research efforts to investigate. Their answers will help resolve the controversy on the oracle property and thus point to the right direction of research on model selection.

%\newpage

\appendix

\newpage

\renewcommand{\theequation}{A.\arabic{equation}}
\setcounter{equation}{0}
\section{Appendix: Proofs of the theorems}
\subsection{Proof of Theorem \ref{asymp_1}}
\begin{m_proof}
Note first that  $\hat{\theta}_{n}\overset
{p}{\rightarrow}\theta$.
For  any $\theta\not =c$, the condition $a_{n}=o(1)$ implies that, for any $\varepsilon>0$, 
\begin{align*}
\Pr_{\theta}(r_{n}\|\breve{\theta}_{n}(c)-\hat{\theta}_{n}\|>\varepsilon)  &
\leq\Pr_{\theta}(\|\hat{\theta}_{n}-c\|\leq a_{n})\leq\Pr_{\theta}(\|\theta-c\|-\|\hat{\theta}_{n}-\theta\|\le a_{n})\\
& =\Pr_{\theta}(\|\hat{\theta}_{n}-\theta\|\geq\|\theta-c\|-a_{n})\rightarrow0\quad
\hbox{as }n\rightarrow\infty.
\end{align*} 
Thus $r_{n}(\breve{\theta}_{n}(c)-\theta)=r_{n}(\breve{\theta}_{n}(c)-\hat{\theta
}_{n})+r_{n}(\hat{\theta}_{n}-\theta)\overset{d}{\rightarrow}Z$.
For $\theta=c$, thanks to $r_{n}a_{n}\rightarrow\infty$,
\begin{equation}\label{Hodges_0}
\Pr_{c}(r_{n}\|\breve{\theta}_{n}(c)-c\|    >\varepsilon)\leq\Pr_{c}%
(\breve{\theta}_{n}\not =c) =\Pr_{\theta_{0}}(r_{n}\|\hat{\theta}_{n}-c\|>r_{n}a_{n}) \rightarrow0.
\end{equation}
This shows $r_{n}(\breve{\theta}_{n}(c)-c)\overset{p}{\rightarrow}0$. 
\end{m_proof}
\subsection{Proof of Theorem
\ref{multi-dimension-super-efficiency}}
\begin{m_proof}
Because we are concerned with the asymptotic distribution of  $\check\theta_{n, b}$ in this section, without loss of generality we can treat the easy case where $V$ is known and $\check\theta_{n, b}$ is defined by 
\begin{equation}\label{part_estimate_2}
\check\theta_{n, b}= \hat\theta_{n,b}+{V}_{bb}^{-1}{V}_{b\bar{b}}(\hat\theta_{n,\bar{b}}-c_{\bar b})
\quad\hbox{and}\quad\check\theta_n(b)=(\check\theta_{n,b}',c_{\bar b}')'
\end{equation}
with the convention $\check\theta_{n,\{1,2,\dots,d\}}=\check\theta_{n}(\{1,2,\dots,d\})=\hat\theta_{n}$.

The first assertion is obvious, so we here only prove (\ref{2.10}). With the definition of $b(\theta)$ in \eqref{b_theta_1}, it is clear that $\theta_{b(\theta)}$  is the sub-vector $(\theta_j: \theta_j\neq c_j)$ of $\theta$ and $\theta_{\bar{b}(\theta)}=(c_j:j\in \bar{b}(\theta))=c_{\bar{b}(\theta)}$. Note that, by \eqref{part_estimate_2}, $\check\theta_{n,b(\theta)}=\hat\theta_{n,b(\theta)}$ and hence $\check\theta_n(b(\theta))=(\check\theta_{n,b(\theta)},c_{\bar b(\theta)})$ are only pseudo-estimators that depend on the unknown parameters $\theta$. However, 
\begin{align}\label{asymptotic_of_check_theta}
r_n(\check\theta_{n,b(\theta)}-\theta_{b(\theta)})
&=V_{b(\theta),b(\theta)}^{-1}(V_{b(\theta)b(\theta)}\quad V_{b(\theta),\bar{b}(\theta)})r_n\left(\begin{array}{c}\hat\theta_{n,{b}(\theta)}-\theta_{b(\theta)}\\
\hat\theta_{n,\bar {b}(\theta)}-\theta_{\bar {b}(\theta)}\end{array}\right)\nonumber\\
&\overset{d}{\rightarrow} V_{{b}(\theta),{b}(\theta)}^{-1}(V_{{b}(\theta){b}(\theta)}\quad V_{{b}(\theta),\bar{b}(\theta)})Z
=\check Z_{b(\theta)},
\end{align}
where the components of $Z$ has been rearranged according to the order in $\theta$ and $\check Z_{b(\theta)}$ is derived from \eqref{X_b} by replacing $b$ with $b(\theta)$. 

For any $\theta$, by comparing \eqref{Hodges' estimate-multi-dimension} and \eqref{part_estimate_2},
\begin{align*}
\Pr_\theta\left(r_n\|\tilde\theta_n(c)-\check\theta_{n}(b(\theta))\|>\varepsilon\right)&=\Pr_\theta\left(r_n\|\tilde\theta_n(c)-(\check\theta_{b(\theta)}', c_{\bar b(\theta)}')'\|>\varepsilon\right)\\
&\leq\Pr_\theta\left(\tilde\theta_n(c)\neq(\check\theta_{b(\theta)}',c_{\bar b(\theta)}')'\right)\\
&\leq\Pr_\theta(b_n(c)\neq b(\theta)).
\end{align*}
Since $\{b_n(c)\neq b(\theta)\}=\bigcup_{j\in b(\theta)}\{|\hat\theta_{nj}-c_j|\leq a_{nj}\}\bigcup_{j\in \bar b(\theta)}\{|\hat\theta_{nj}-c_j|> a_{nj}\}$,
\begin{align*}
\Pr_\theta&\left(r_n\|\tilde\theta_n(c)-\check\theta_{n}(b(\theta))\|>\varepsilon\right)\\
&\leq\Pr_\theta\left(\bigcup_{j\in b(\theta)}\{|\hat\theta_{nj}-c_j|\leq a_{nj}\}\bigcup_{j\in \bar b(\theta)}\{|\hat\theta_{nj}-c_j|> a_{nj}\} \right)\\
&\leq\sum_{j\in b(\theta)}\Pr_\theta\{|\hat\theta_{nj}-c_j|\leq a_{nj}\}+\sum_{j\in \bar b(\theta)}
\Pr_\theta\{|\hat\theta_{nj}-c_j|> a_{nj}\} \\
&\leq \sum_{j\in b(\theta)}\Pr_\theta\{|\hat\theta_{nj}-\theta_j|\geq |\theta_j-c_j|-a_{nj}\}
+\sum_{j\in \bar b(\theta)}\Pr_\theta\{r_n|\hat\theta_{nj}-c_j|> r_na_{nj}\}\\
& \rightarrow0 \quad\hbox{as }n\rightarrow\infty 
\end{align*}
under the conditions on $a_n$. Thus 
$r_n(\tilde\theta_n(c)-\check\theta_{n}(b(\theta))=o_p(1)$. Combining this with \eqref{asymptotic_of_check_theta}, we get
\begin{equation}\label{asymp_2}
r_n(\tilde\theta_n(c)-\theta)=r_n(\tilde\theta_n(c)-\check\theta_{n}(b(\theta)))
+r_n(\check\theta_{n}(b(\theta))-\theta)\overset{d}{\rightarrow} \left(\begin{array}{c}\check Z_{b(\theta)}\\
0\end{array}\right)
\end{equation}
under $\Pr_\theta$. Next examine the case $b(\theta)=\emptyset$. Analogous to \eqref{Hodges_0}, 
under $\Pr_c$ and the condition 
$r_n\min_{1\leq j\leq d}\limits a_{nj}\rightarrow \infty$, 
\begin{equation}\label{Hodges_1}
\Pr_{c}(r_{n}\|\tilde{\theta}_{n}(c)-c\|    >\varepsilon)\leq\Pr_{c}%
(\tilde{\theta}_{n}\not =c) \leq\sum_{j=1}^d\Pr_{c}(r_{n}|\hat{\theta}_{nj}-c|>r_{n}a_{nj}) \rightarrow0
\quad\hbox{as }n\rightarrow\infty.
\end{equation}
Hence $r_{n}|\tilde{\theta}_{n}(c)-c| \overset {d}{\rightarrow}0$. The proof is then complete.
\end{m_proof}
\end{document}